\makeatletter

\def\eqalign#1{\,\vcenter{\openup\jot\m@th
  \ialign{\strut\hfil$\displaystyle{##}$&$\displaystyle{{}##}$\hfil
      \crcr#1\crcr}}\,}
\def\eqalignno#1{\displ@y \tabskip\@centering
  \halign to\displaywidth{\hfil$\displaystyle{##}$\tabskip\z@skip
    &$\displaystyle{{}##}$\hfil\tabskip\@centering
    &\llap{$##$}\tabskip\z@skip\crcr
    #1\crcr}}
\def\leqalignno#1{\displ@y \tabskip\@centering
  \halign to\displaywidth{\hfil$\displaystyle{##}$\tabskip\z@skip
    &$\displaystyle{{}##}$\hfil\tabskip\@centering
    &\kern-\displaywidth\rlap{$##$}\tabskip\displaywidth\crcr
    #1\crcr}}

\makeatother

\newdimen\pixel \pixel=.00333333 true in

\documentclass[11pt]{article}
\usepackage{epsfig}
\usepackage{fullpage}

\makeatletter
\def\bigpar{\bigbreak\@afterindentfalse\@afterheading\ignorespaces}
\def\medpar{\medbreak\@afterindentfalse\@afterheading\ignorespaces}
\def\smallpar{\smallbreak\@afterindentfalse\@afterheading\ignorespaces}

\newlength{\saveindent}
\newenvironment{proof}%
      {\bigpar{\bf Proof:}\ 
             \setlength{\saveindent}{\parindent} 
                       \ignorespaces}%
      {\stopproof\ignorespaces\bigbreak \setlength{\parindent}{\saveindent}}

      {\bigpar{\bf Proof:}%
             \setlength{\saveindent}{\parindent} 
                       \,\ignorespaces}%
      {\ignorespaces\bigbreak \setlength{\parindent}{\saveindent}}

      {\bigpar{\bf Proof:}\ %
             \setlength{\saveindent}{\parindent} 
                       \ignorespaces}%
      {\ignorespaces\bigbreak \setlength{\parindent}{\saveindent}}

      {\bigpar{\bf#1:}\ %
             \setlength{\saveindent}{\parindent} 
                       \ignorespaces}%
      {\stopproof\ignorespaces\bigbreak \setlength{\parindent}{\saveindent}}

      {\smallpar{\bf Remark:}\ 
                       \ignorespaces}%
      {\stopproof\ignorespaces\medbreak \setlength{\parindent}{\saveindent}}

\newenvironment{remark*}%
      {\smallpar{\bf Remark:}\ 
                       \ignorespaces}%
      {\ignorespaces\medbreak \setlength{\parindent}{\saveindent}}

      {\smallpar{\bf Remarks:}\ 
                       \ignorespaces}%
      {\stopproof\ignorespaces\medbreak \setlength{\parindent}{\saveindent}}

\newenvironment{remarks*}%
      {\smallpar{\bf Remarks:}\ 
                       \ignorespaces}%
      {\ignorespaces\medbreak \setlength{\parindent}{\saveindent}}

      {\smallpar{\bf Remark:}\ %
                       \ignorespaces}%
      {\stopproof\ignorespaces\medbreak \setlength{\parindent}{\saveindent}}

      {\smallpar{\bf Remarks:}\ %
                       \ignorespaces}%
      {\stopproof\ignorespaces\medbreak \setlength{\parindent}{\saveindent}}
\makeatother

\newtheorem{theorem}{Theorem}
\newtheorem{lemma}[theorem]{Lemma}

\newtheorem{example}{Example}

\def\begex{\begin{example}\parindent=0pt \rm}
\def\endex{\end{example}}
\def\square{\vbox{\hrule height.2pt\hbox{\vrule width.2pt height5pt \kern5pt
                                   \vrule width.2pt} \hrule height.2pt}}
\def\stopproof{\hfill \square \smallskip}

\def \is {{i^*}}

\def \ltwo#1{{||#1||_2}}
\def \ss#1{{||#1||^2_2}}
\def\sfrac#1#2{{\textstyle{#1 \over #2}}}

\def\quarter{{\textstyle{1\over4}}}

\def \ps {{ \sum_x p_x^2}}

\def \r {{\bf R}}
\def \var {{ \rm var }}
\def \r {{ \cal F}}

\def\r|{{\Bigr\vert}}
\def\l|{{\Bigl\vert}}

\def \lov {{Lov\'asz }}
\def\phi {\Phi}

\def\e{\epsilon}

\def \cov {{\rm cov}}

\def\tmix{\tau_{\rm mix}}
\def\tmixp{\tau'_{\rm mix}}

\def\varepsilon{\mathchar"122 }

\def \one {{\mathbf 1}}
\def \chi {{\mathbf 1}}

\def\p {{ \bf P}}
\def\P {{ \bf P}}

\def\e {{ \bf {E}}}

\def \ns {{\tt nextstate()}}

\def\Sq{{\cal S}_q}
\def\Sq-{{\cal S}_{q-1}}

\def\u{{\cal U   }}

\def\one{{\mathbf 1}}

\newcommand{\be}{\begin{equation}}
\newcommand{\ee}{\end{equation}}
\newcommand{\lab}{\label}
\begin{document}
\title{The birthday problem and Markov chain Monte Carlo}
\author{
  {\sc Itai Benjamini}\thanks{Weizmann Institute,
\newline Email:
{\tt itai.benjamini@weizmann.ac.il}}
and
{\sc Ben Morris}\thanks{University of California, Davis
\newline Email:
{\tt morris@math.ucdavis.edu}.} }
\maketitle
\begin{abstract}
We study the problem of generating a sample from the stationary
distribution of a Markov chain, given a method to simulate the
chain. We give an approximation algorithm for the case of a random
walk on a regular graph with $n$ vertices that runs in expected
time $O^*(\sqrt{n} \tmix)$, where $\tmix$ is the mixing time of
the chain in $L^2$. This is close to the best possible,
since
$\sqrt{n}$ is a lower bound on the worst-case expected
running time of any algorithm.
\end{abstract}
\section{Introduction} \lab{intro}
Often the only feasible method for sampling from a complex
distribution is to simulate a suitably chosen Markov chain for
sufficiently many steps. Procedures based on this idea, called
Markov chain Monte Carlo, have been applied to a number of
problems such as approximating the  permanent \cite{jsv},
computing volumes \cite{dfk} and integrals \cite{integ}, and
approximate counting \cite{v}. In order to find out how many steps
one needs to simulate the chain, it is necessary to determine the
{\it mixing time}, i.e., the number of steps necessary to bring
the distribution close to stationary. This analysis is often
complex and has to be tailored to the specific type of Markov
chain under consideration.

Another, related line of research has been pursued (see \cite{pw}
for background and see also \cite{agt, aldous, lw}): the Markov
chain is generic, its transition probabilities are not known, and
the algorithm is given a procedure to generate the next state of
the chain based on the current state. Since the chain is
arbitrary, the algorithm cannot have any advanced knowledge of the
mixing time. Aldous \cite{aldous}, comes by time $O(\tau
\mbox{Poly} (1/\epsilon))$ to within $\epsilon$ from the
stationary distribution in total variation.
This was improved by
the ``cycle popping'' algorithm invented by Propp
and Wilson \cite{pw}, which  runs in expected time $O(\tau)$, where
$\tau$ is the mean hitting time (expected number of steps the
chain takes to get from $X$ to $Y$ if $X$ and $Y$ are chosen
independently according to the stationary distribution). In
\cite{lw} a random stopping rule for exact sampling from an
unknown Markov chain is given where the expected number of steps is
$6 \tau^4$. It is easy to see that any algorithm for an arbitrary
$n$-state chain must have running time at least $n$: the algorithm
must visit every state, since an unvisited state $x$ could
potentially have a holding probability $p(x,x)$ very close to $1$
and hence a very high stationary probability.  In this note, we
show that if the chain is a random walk on a regular graph then
there is an approximation algorithm whose running time can be much
less. Our algorithm runs in expected time $O^*(\sqrt{n} \tmix)$,
where $\tmix$ is the $L^2$ mixing time of the chain; this can be
much smaller than $n$.
We note that
our algorithm produces a sample with some error (i.e., the
distribution is not exactly stationary, although it can be made
arbitrarily close), whereas the algorithms described in \cite{pw,
lw, agt} produce an exact sample. In a somewhat similar spirit of
observing a random walk on an unknown graph, \cite{bglrt} studies
what can be learned by knowing polynomially many return times to a
fixed vertex of a simple random walk on a regular graph.

 Our analysis is based on a variation of the standard {\it birthday
 problem}. Roughly speaking,
in a world where there are $n$ possible birthdays, the number of
 people that you need to pick to be confident that at least two have the
 same birthday is of the order $\sqrt{n}$. It turns out that
in a similar vein,
if order
$\sqrt{n}$ copies of a Markov chain with uniformly stationary distribution
are run for much less than the mixing time then there is likely to be
 a match, whereas if they are run for much more than the mixing time
 there is a good chance for no match.
(See Section \ref{main} for a precise formulation of this.) This
forms the basis for our algorithm. A similar idea was used by
Goldreich and Ron \cite{gr}, as a suggestion for a possible
sublinear tester for expansion.

\section{Results}
\label{main}
\subsection{The Problem}
The algorithm is given an irreducible, aperiodic
$n$ state Markov chain as input.
More precisely, the
algorithm is given the number of states $n$, a starting state $x_0$,  and
a procedure \ns, which, given a state $x$ of the chain outputs
state $y$ which is one step of the chain starting from $x$.
The problem is to
generate a random state according to the stationary
distribution.

We will aim for an $\epsilon$-approximation algorithm, that is, an
algorithm that generates a random state within total variation
distance $\epsilon$ of the stationary distribution.

\subsection{Main Theorem}
Let $p^m(\,\cdot \,,\, \cdot )$
be transition probabilities for an irreducible, aperiodic,
doubly stochastic $n$-state Markov chain on state space $V$.
Let $\u$ denote the uniform distribution over $V$ and
for functions $\mu$ on $V$, let
$\ltwo{\mu} = \Bigl(\sum_{x \in V} \sfrac{1}{n}\mu(x)^2\Bigr)^{1/2}$ denote
the norm of $\mu$ in $L^2(\u)$.
For $\epsilon > 0$, let $\tau(\epsilon) =
\min\{n: \ss{n p^n(x_0, \, \cdot) - 1} \leq \epsilon\}$.
Denote the mixing time in $L^2$ by
$$
\tmix  = \max_x \min\{n: \ss{n p^n(x, \, \cdot) - 1} \leq 1/e\}.
$$
Our main result is
the following theorem.

\begin{theorem}
Suppose that the Markov chain is a random walk on a regular, connected graph
with degree at most $n$. Then there is an algorithm that returns
a sample within total variation distance $\epsilon$
in expected time $O( \sqrt{n}
\log(n/\epsilon) \tmix)$.
\end{theorem}
\medskip

\noindent {\bf Remark:} To obtain a {\it lower bound} for the
running time, we can consider random walk on the complete graph.
If an algorithm simulates the chain for less than order $\sqrt{n}$ steps, it
is likely to see only distinct states of the chain, hence it
couldn't ``tell the difference'' between
the chain and random walk on
two
complete graphs of size $n/2$ joined by a single edge.
(Note that in the second case
one would have to do at least order $n$ simulated steps
to get an almost uniform sample.)
So any algorithm needs at least order $\sqrt{n}$
steps, which is order $\sqrt{n} \tmix$
steps since the mixing time is a constant.
We believe that for analogous reasons this would still
hold
(i.e., any algorithm would need order
$\sqrt{n} \tmix$ steps)
when the complete graph is
replaced by a random $d$-regular graph,
but it seems harder to prove this.
\medskip

\begin{proof}

\noindent The algorithm is as follows. Define $A_n = n^4
\log(2n/\epsilon)$ and $l = \lceil 1 + {512 \sqrt{n} \over
\epsilon^2} \rceil$. Iterate the following procedure for $i \in
\{1, 2, \dots\}$ until stopped.

Let $m = 8\log(2A_n/\epsilon)$, and perform the following experiment $m$
times. Simulate $l$ copies of the Markov chain starting at $x_0$
for $2^i$ steps,
generating $l$ samples $X_1, \dots, X_l$.
To avoid the possibility that the Markov chain has an eigenvalue close
to
$-1$, we implement a holding probability of $1/n$ to each state; i.e.,
each step we do nothing with probability $1/n$, otherwise
simulate a step of the chain.
This doesn't change the stationary distribution and ensures that all
eigenvalues are at least $-1 + \sfrac{2}{n}$.
Let $\delta = \epsilon^2$.
Let
\begin{equation}
\label{zdef}
Z = \sum_{k < m}
\one(X_k = X_m),
\end{equation}
and if $Z \le (1 + {\delta/2}) {l \choose 2} {1 \over n}$
then count the experiment as a success; otherwise count it as a
failure. If at least $m/2$ of the experiments are successful,
or if $i = \lceil \log_2 A_n \rceil$
then stop; otherwise continue with the next value of $i$.

Let $i'$ be the value of $i$ when the above procedure terminates.
We claim that with high probability after $2^{i'}$  steps the
chain is very mixed. Hence the algorithm can run another
independent simulation of the chain for $2^{i'}$ steps and the
result is an almost uniform sample from the state space.

More precisely, let $\mu(i) = p^{2^i}(x_0, \cdot)$.
We will show that with probability at least $1 - \epsilon/2$,
the value of $i'$ is large enough so that
$\ss{n\mu(i') - 1} \leq \delta$.
\bigskip

\noindent {\bf Analysis of the algorithm.}  By a {\it conductance}
bound (see, e.g., \cite{js}), the spectral gap for the chain must
be at least $1/n^4,$ and hence $\tau(\delta/2) \leq n^4
\log(2n/\delta) = A_n$. Thus if $i' = \lceil \log_2 A_n \rceil$
then $2^{i'} \geq A_n$, which implies that $\ss{n\mu(i') - 1 }
\leq \delta/2$.

Next we have to bound the probability that the algorithm stops
early on a value of $i$ such that $\ss{ n\mu(i) - 1} \geq
\delta/2$. Fix $i < \lceil \log_2 A_n \rceil$ and let $p_x$ be
the probability that the chain is at $x$ after $2^i$ steps.
Cauchy-Schwarz gives $\sum_x p_x^2 \geq 1/n$. It follows that if
$Z$ is defined as in (\ref{zdef}) then $\e(Z) = {l \choose 2} \ps
\geq {1 \over n}{(l-1)^2 \over 2},$ and hence \be \lab{b1} {1
\over \e(Z)} \leq {\delta^2 \over 512}. \ee We also have \be
\lab{b2} {2 \sqrt{n} \over l} \leq {\delta^2 \over 512}. \ee
Note that $\e(Z) = {l \choose 2} {1 \over n}(1 + \ss{n \mu(i) - 1})$.
Thus
if $\ss{ n\mu(i) - 1} \geq \delta/2$, then
$\e(Z) \geq {l \choose 2} {1 \over n}
(1 + \delta/2)$
and hence
\begin{eqnarray}
\P\Bigl( Z \leq {l \choose 2} {1 \over n} (1 + \delta/2)\Bigr)
&\leq& \P\Bigl( Z \leq \e(Z) \Bigl( {1 + {\delta}/{2} \over 1 +
\delta}
\Bigr)\Bigr) \\
&\leq& \P\Bigl( Z \leq \e(Z) \Bigl( {1 - {\delta}/{4}}
\Bigr)\Bigr) \\
&\leq& {\var(Z) \over \Bigl[ \sfrac{\delta^2}{16} \e(Z)^2     \Bigr]} \\
\label{bound2}
&\leq& {16}{\delta^{-2}} \cdot {1 + {2 \sqrt{n} \over l}\e(Z)
\over \e(Z)},
\end{eqnarray}
where the second inequality follows from the fact that
${1 + \delta/2 \over 1 + \delta} \leq 1 - \delta/4$ whenever
$\delta \leq 1$, the third line is Chebyshev's inequality,
and the fourth line is Lemma \ref{varlemma} from the Appendix.
The upper bound (\ref{bound2}), and hence the probability of
success is at most
\[
{16 \over \delta^2}
\Bigl({1 \over \e(Z)} + {2 \sqrt{n} \over l}\Bigr)
\leq {1 \over 4}.
\]
Thus Hoeffding's bounds imply that the probability of at least
$m/2$ successes in stage $i$ is at most
$e^{-m/8} \leq \epsilon/2A_n$. Since the number of stages can never
be more than $A_n$, the probability that $i'$ is such that
$\ss{n\mu(i') - 1} \geq \delta$ is at most $\epsilon/2$.
Recall that $\delta = \epsilon^2$ and note that
the total variation distance
$\Vert \mu(i') - \u \Vert_{TV}
= 2 || n\mu(i') - 1 ||_1$.
Hence, when $\ss{n\mu(i') - 1} \leq \delta$ we
have $|| n\mu(i') - 1 ||_1 \leq
\ltwo{ n\mu(i') - 1}
\leq \delta$ and hence
the total variation distance
$\Vert \mu(i') - \u \Vert_{TV} \leq \epsilon/2$.
It
follows that the algorithm will generate a random sample
within total variation
distance $\epsilon$ of uniform.\\

\noindent {\bf Running time.}
Define $\tmixp = \tau(\delta/4)$ and let
$\is = \min\{i: 2^i \geq \tmixp\}$.
If $i \leq \is$ then the
number of steps in stage $i$ is at most $C(l 2^i \log(2A_n/\epsilon)$
for a universal constant $C$.
(We assume that the values of the chain are given as
0-1 strings, whose lengths we treat as constant.
If we  store the values of $X_1, \dots, X_l$ in a binary tree,
then we can
count the number of matches among them in $O(l)$ time.)
Summing this over $i \leq \is$ shows that the number of steps
corresponding to $i \leq \is$ is
$O(l 2^{\is} \log(2A_n/\epsilon)) =
O(l \tmixp \log(2A_n/\epsilon))$. Suppose that $i \geq \is$.
Then $2^i \geq \tmixp$ and hence $\ss{\mu(i) - \u}  \leq
\delta/4$.
Step $i+1$ occurs only if
there are more than $m/2$ failures in step $i$. Note that
$\e(Z) \leq {l \choose 2} {1 \over n}
(1 + {\delta}/{4})$
and hence
\begin{eqnarray}
\P\Bigl( Z > {l \choose 2} {1 \over n} (1 + {\delta}/{2})\Bigr)
&\leq& \P\Bigl( Z > \e(Z) \Bigl( {1 + {\delta}/{2} \over 1 +
{\delta}/{4}}
\Bigr)\Bigr) \\
&\leq& \P\Bigl( Z > \e(Z) \Bigl( {1 + \sfrac{\delta}{8}}
\Bigr)\Bigr) \\
&\leq& {\var(Z) \over \Bigl[ \sfrac{\delta^2}{64} \e(Z)^2     \Bigr]} \\
&\leq&
\label{bound1}
\sfrac{64}{\delta^2} \cdot {1 + {2 \sqrt{n} \over l}\e(Z)
\over \e(Z)} \leq \quarter,
\end{eqnarray}
\\
where the third line is Chebyshev's inequality
and the second line uses the fact that
${1 + \delta/2 \over 1 + \delta/4} \geq 1 + \delta/8$ whenever
$\delta \leq 1$.
Thus Hoeffding's bounds give
$\P(\mbox{stage $i + 1$ occurs}) \leq
\P(\mbox{more than $m/2$ failures}) \leq
e^{-m/8} \leq \epsilon/2A_n$. Since the maximum number of steps in
any stage is $O(A_n l \log(2A_n/\epsilon))$, summing the above bound
over $i$ gives an $O(l \log(2A_n/\epsilon))$
bound. Adding everything up gives a total expected running time
of $O(l \tmixp \log(2A_n/\epsilon))
=
O(l \tmix \log(n/\epsilon))$.
\bigskip
\end{proof}
\bigskip

\noindent {\bf Remarks:}

\begin{itemize}
\item
By running $\sqrt{n}$ samples after the mixing time is estimated
instead of just one, the algorithm could actually
produce $\sqrt{n}$ samples and the expected
running time would still be
$O^*(\sqrt{n} \tmix)$.

\item
The assumption that the degree is at most $n$ can be relaxed;
it is only used to get a
poly($n$) upper bound for the mixing time (in order to bound $A_n$).

\item
In \cite{abls} it was observed that a non-backtracking random walk
mixes (up to a factor two) faster. Thus, in the setting where
the algorithm can determine the set of neighbors of a
state $x$, one can very slightly
reduce the randomness used, as well as the running time, by
replacing the simple random walk with a non-backtracking random
walk. Also simulating the $l$ copies of the Markov chain at each
stage, can of course be done in parallel.
\end{itemize}

\section{Appendix}
The following bound on the variance of $Z$ was needed.
\begin{lemma}
\label{varlemma}
Let $X$ be a random variable taking values in $\{1, \dots, n\}$,
and let $p_i = \p(X = i)$ for $1 \leq i \leq n$.
Let
$X_1, \dots, X_l$ be independent copies of $X$ and
let $Z = \sum_{i < j} \one(X_i = X_j)$.
Then
\begin{enumerate}
\item
$\e(Z) = {l \choose 2} \sum_i p_i^2$.
\item
$\var(Z)
\leq \e(Z) \Bigl( 1 + {2 \sqrt{n} \over l} \e(Z)\Bigr).$
\end{enumerate}
\end{lemma}
\begin{proof}
Part 1 is obvious. For part 2, let $M = \max_i p_i$. Clearly,
$\sum_i p_i^2 \geq M^2$, and Cauchy-Schwarz gives $\sum_i p_i^2 \geq
{1 \over n}$. Hence
\[
\sum_i p_i^2 \geq \sqrt{ M^2 \cdot {1 \over n} } = {M \over \sqrt{n}}.
\]
It follows that
\be
\lab{bound}
\sum_i p_i^3 \leq M \sum_i p_i^2 \leq \sqrt{n} \Bigl( \sum_i p_i^2
\Bigr)^2.
\ee
\\
\\
For $1 \leq i<j \leq n$, let $Z_{i,j} = \one(X_i = X_j)$, so that
$Z = \sum_{i < j} Z_{i,j}$.
Note that $\cov(Z_{i,j}, Z_{k,l}) = 0$ if $i,j,k,l$ are distinct.
Thus,
\begin{eqnarray*}
\var(Z) &=& \cov\Bigl( \sum_{i \neq
j} Z_{i,j}, \sum_{i \neq j} Z_{i,j}\Bigr) \\
&=& {l \choose 2} \cov\Bigl(Z_{1,2}, Z_{1,2} \Bigr)
+ 2 {l \choose 3} \cov(Z_{1,2},Z_{1,3}) \\
&\leq& {l \choose 2} \e(Z_{1,2}) + 2 {l \choose 3} \e(Z_{1,2}Z_{1,3}) \\
&\leq& \e(Z) + 2{l \choose 3} \sum_i p_i^3.
\end{eqnarray*}
Combining this with equation (\ref{bound})
and part 1 of the lemma yields
\begin{eqnarray*}
\var(Z)
&\leq& \e(Z) + {l \choose 3} \sqrt{n} \Bigl( \sum_i p_i^2 \Bigr)^2\\
&=& \e(Z) \Bigl( 1 + {2 {l \choose 3} \sqrt{n} \over
{l \choose 2}^2} \e(Z)
\Bigr) \\
&\leq& \e(Z) \Bigl( 1 + {2 \sqrt{n} \over l } \e(Z)\Bigr), \\
\end{eqnarray*}
establishing
part 2 of the lemma.
\end{proof}
\medskip

\noindent {\bf Acknowledgement:} Thanks to Noam Berger for a
useful discussion.


\begin{thebibliography}{99}
\bibitem{aldous}
Aldous, D. On simulating a Markov chain stationary distribution
when the transition probabilities are unknown.  {\it Discrete
Probability and Algorithms,}   {\it IMA Volumes in Mathematics and
its Applications,} {\bf 72} Springer-Verlag, (1995), pp.~1--9.

\bibitem{abls}
Alon, N. Benjamini, I. Lubetzky, E. and S. Sodin, Non-backtraking
random walk mixes faster. Preprint (2006).


\bibitem{agt}
Asmussen, S., Glynn, P., and Thorisson, H. Stationary detection in
the initial transient problem. {\it ACM Transactions on Modeling
and Computer Simulation.} {\bf 2} (1992), pp.~130--157.


\bibitem{bglrt}
Benjamini, I. Kozma, G. \lov, L. Romik, D. and Tardos, G. Waiting
for a Bat to Fly By (in Polynomial Time). {\it Combinatorics,
Probability and Computing} {\bf 15} (2006),  pp.~673-683.

\bibitem{gr}
Goldreich, O., and Ron, D. On testing expansion in bounded-degree
graphs.  {\it ECCC}, TR00-020, (2000).


\bibitem{js} Jerrum, M.~and Sinclair, A.
Approximating the permanent. {\it SIAM Journal on Computing.} {\bf
18} (1989), pp.~1149--1178.

\bibitem{jsv} Jerrum, M., Sinclair, A., and Vigoda, E.
A polynomial-time approximation algorithms for the permanent of a
matrix with non-negative entries. {\it Journal of the ACM.} {\bf
51} (2004), pp.~671--697.

\bibitem{dfk} Dyer, M, Frieze, A., and Kannan, R.
A random polynomial time algorithm for approximating the volume of
convex sets. {\it Journal of the ACM.} {\bf 38} (1991), pp.~1--17.

\bibitem{integ}
Frieze, A., Kannan, R., and Polson, N. Sampling from log-concave
distributions. {\it Annals of Applied Probability.} {\bf 4}
(1994), pp.~812--837.

\bibitem{v}
Jerrum, M., Valiant, L., and Vazirani, V. Random generation of
combinatorial structures from a uniform distribution. {\it
Theoretical Computer Science.} {\bf 43} (1986), pp.169--188.

\bibitem{lw}
\lov, L.~and Winkler, P.
Exact mixing in an unknown Markov chain.
{\it Electronic Journal of Combinatorics.}
 {\bf 2} (1995).
Paper \#R15.

\bibitem{pw}
Propp, J.~and Wilson, D. How to get a perfectly random sample from
a generic Markov chain and generate a random spanning tree of a
directed graph. {\it Journal of Algorithms.} {\bf 27} (1998),
pp.170--217.

\end{thebibliography}
\end{document}